\documentclass[11 pt, oneside, reqno]{amsart}
\usepackage[top=1.5in,bottom=1in,left=1in,right=1in]{geometry}
\usepackage{amsmath,amssymb,amsthm,amsopn,extarrows,accents,faktor,enumitem,stmaryrd,mathtools}  
\usepackage[linktocpage=true]{hyperref}
\hypersetup{
    colorlinks = true,
    allcolors = {blue}}
\usepackage[dvipsnames]{xcolor}

\usepackage{setspace}
\usepackage{tikz-cd,tikz}
\usetikzlibrary{matrix,arrows,decorations.pathmorphing}

\newtheorem{theorem}[equation]{Theorem}
\newtheorem*{theorem*}{Theorem}
\newtheorem{proposition}[equation]{Proposition}

\theoremstyle{definition}
\newtheorem{definition}[equation]{Definition}
\newtheorem*{remark}{Remark}
\newtheorem{example}[equation]{Example}
\newtheorem{exercise}{Exercise}

\numberwithin{equation}{section}

\newcommand{\hooklongrightarrow}{\lhook\joinrel\longrightarrow}

\newcommand{\g}{\mathfrak{g}}
\renewcommand{\b}{\mathfrak{b}}
\newcommand{\h}{\mathfrak{h}}
\renewcommand{\t}{\mathfrak{t}}
\newcommand{\p}{\mathfrak{p}}
\renewcommand{\l}{\mathfrak{l}}
\renewcommand{\d}{\mathfrak{d}}
\renewcommand{\a}{\mathfrak{a}}
\newcommand{\n}{\mathfrak{n}}
\renewcommand{\k}{\mathfrak{k}}
\newcommand{\uu}{\mathfrak{u}}
\DeclareMathOperator{\Ad}{Ad}
\DeclareMathOperator{\ad}{ad}

\DeclareMathOperator{\id}{Id}

\newcommand{\xo}{\mathring{X}}
\renewcommand{\O}{\mathcal{O}}
\newcommand{\T}{\mathcal{T}}
\renewcommand{\L}{\mathcal{L}}
\newcommand{\Gs}{\widetilde{G}}
\newcommand{\Ts}{\widetilde{T}}
\newcommand{\Bs}{\widetilde{B}}
\newcommand{\Gbar}{\overline{G}}
\newcommand{\mubar}{\overline{\mu}}
\DeclareMathOperator{\stab}{Stab}
\DeclareMathOperator{\gr}{Gr}

\makeatletter
\def\l@subsection{\@tocline{2}{0pt}{3pc}{5pc}{}}
\makeatother

\title{Wonderful varieties with a view towards Poisson geometry}
\author{Ana B\u{a}libanu}
\address{Department of Mathematics, Harvard University, 1 Oxford Street, Cambridge, MA  02138, USA}
\email{ana@math.harvard.edu}
\date{}

\begin{document}
\maketitle

\begin{abstract}
These are expanded lecture notes from the author's minicourse at the 2022 Poisson Geometry Summer School, which took place at the Centre de Recerca Matem\`atica in Barcelona, Spain. After giving a general introduction to wonderful varieties, and an explicit construction of the wonderful compactification of a semisimple adjoint group, we outline several connections to Poisson geometry and to varieties of Lagrangian subalgebras. This survey is intended to be accessible to readers without an extensive background in algebraic geometry.
\end{abstract}

\tableofcontents

A homogeneous space is an algebraic variety $\xo$ on which an algebraic group $G$ acts transitively. To study $\xo$, it is often useful to study its equivariant compactifications---that is, those compactifications in which the action of $G$ extends to the boundary. Such compactifications model the equivariant behavior of $\xo$ ``at infinity,'' where many interesting features of its geometry become apparent. For several large classes of homogeneous spaces, such compactifications can be described entirely by discrete combinatorial data. The most famous example of such a classification is the theory of toric varieties, in which equivariant compactifications of a torus are indexed by fans or polytopes. 

The goal of these notes is to give an introduction to an important class of equivariant compactifications of homogeneous spaces which are called wonderful varieties. We begin with an overview that places these objects into the broader setting of spherical varieties, and explains some of the motivation behind the development of this theory. Then we construct the wonderful compactification of a semisimple algebraic group, and we describe some of its remarkable structure. In the final three sections, we connect these constructions to varieties of Lagrangian subalgebras and to Poisson geometry.

\subsection*{A note on background} A main goal of the minicourse on which these notes are based was to make this classical topic accessible to an audience without a broad knowledge of algebraic geometry. These notes assume only some basic notions of algebraic groups, Lie algebras, group actions, and algebraic varieties, and provide a long list of references for detailed introductions of these subjects. We have also made an effort to give concrete examples and intuitive heuristics to illustrate the geometry of the constructions that follow. In particular, readers who are not very familiar with algebraic varieties will lose nothing by thinking of algebraic groups as complex Lie groups, and of smooth algebraic varieties as holomorphic manifolds. 

\subsection*{Acknowledgements} I would like to thank the organizers and the scientific committee of Poisson 2022 for the invitation to give the minicourse on which these notes are based. I am also grateful to Yu Li for his help in running this minicourse and for many interesting discussions.\\

%
%
%
%
%
%
%
\section{Introduction}
Let $H$ be a connected algebraic group over the complex numbers. A \emph{$H$-homogeneous space} is an algebraic variety $\xo$ equipped with a transitive action of $H$. By fixing a point $x\in\xo$ we obtain an identification
\[\xo\cong H/K,\]
where $K=\stab_H(x)$ is the stabilizer of $x$ in $H$.

\begin{definition}
\label{def:equiv}
An \emph{equivariant embedding} of the homogeneous space $\xo$ is a normal $H$-variety $X$ together with an $H$-equivariant embedding 
\[\xo\hooklongrightarrow X\]
whose image is open and dense. An \emph{equivariant compactification} of $\xo$ is an equivariant embedding which is a projective variety.
\end{definition}

\noindent For details on the study of equivariant embeddings, we refer the reader to the extensive introduction \cite{tim:11} by Timashev.

While equivariant embeddings may fail to be smooth, in this exposition we will require them to be \emph{normal}. This means that the local ring at every point is an integrally closed domain. One consequence of normality is that algebraic local functions on the smooth locus of $X$ extend to $X$. In particular, as in Hartog's theorem this implies that $X$ has no singularities in codimension one.

Given a homogeneous space, one can ask for a classification of all its equivariant embeddings. The most well-known answer to this type of question is the theory of \emph{toric varieties}, which are equivariant embeddings of the complex torus
\[H=(\mathbb{C}^*)^n.\]
These are parametrized by discrete combinatorial objects known as \emph{fans}, which are collections of convex rational polyhedral cones in $\mathbb{R}^n$ closed under taking intersections and faces. All the geometric features of toric varieties, like their orbit structure, their cohomology, the morphisms between them, etc., are indexed by combinatorial data encoded in their fans. For more background, we refer to the classical introduction \cite{ful:93} by Fulton.

For more general choices of $H$ the collection of all equivariant embeddings of a given homogeneous space is often intractably large. For example, even the equivariant compactifications of the apparently simple additive group $\mathbb{C}^n$ are not ``discrete''---instead, they are parametrized by moduli spaces \cite{has.tsi:99}. However, a good theory of equivariant embeddings does exist for certain relatively large classes of homogeneous spaces.

The first distinguished class that we will be interested in consists of symmetric spaces. To define them, we will need the notion of an \emph{involution} of $H$, which is a nontrivial algebraic group automorphism $\sigma:H\longrightarrow H$ that satisfies $\sigma\circ\sigma=\id_H.$

\begin{definition}
\label{symspaces}
A homogeneous space $H/K$ is a \emph{symmetric space} if $K$ is an open algebraic subgroup of the fixed-point subgroup
\[H^\sigma\coloneqq\{h\in H\mid \sigma(h)=h\}\]
of an involution $\sigma$. Since algebraic subgroups are always closed, this condition is equivalent to requiring that
\[(H^\sigma)^\circ\subset K\subset H^\sigma.\]
\end{definition}

\begin{remark}
If $H$ was a compact real Lie group, Definition \ref{symspaces} would coincide with the usual ``Riemannian'' definition of a symmetric space, as a Riemannian manifold which has an inversion symmetry centered at every point.
\end{remark}

\begin{example}
Suppose that $H=G\times G$ is a direct product of two copies of a group $G$, and consider the involution
\begin{align*}
\sigma:G\times G&\longrightarrow G\times G\\
		(g,h)&\longmapsto (h,g).
		\end{align*}
Its locus of fixed points is 
\[H^\sigma=G_\Delta\coloneqq\{(g,g)\in G\times G\mid g\in G\},\]
which is the stabilizer of the identity element under the transitive action of $G\times G$ on $G$ given by
\[(g,h)\cdot a= gah^{-1}.\]
Therefore we get an identification
\[G\cong (G\times G)/G_\Delta\]
which shows that the group $G$ itself is a symmetric $G\times G$-space.
\end{example}	

\begin{exercise}
\label{ex1}
Suppose that $U\subset H$ is a connected normal unipotent subgroup of $H$.
\begin{enumerate}[label=(\alph*)]
\item Let $\sigma$ be an involution of $H$. Show that $U\subset (H^\sigma)^\circ.$
\item Conclude that, if $\xo$ is a symmetric space for the action of $H$, then $\xo$ is also a symmetric space for the action of a quotient of $H$ which has no nontrivial connected normal unipotent subgroups.
\end{enumerate}
\end{exercise}

Groups that have no nontrivial connected normal unipotent subgroups are called \emph{reductive} groups. One can show that a group $H$ is reductive if and only if all its finite-dimensional representations are completely reducible. Complex tori, the general linear group $GL_n$, and the complex algebraic groups $SL_n$, $SO_n$, and $Sp_{2n}$ are all examples of reductive groups, and we refer to \cite{hum:75} for more background on their structure theory. A special class of reductive groups consists of the \emph{semisimple} groups, which are groups with no nontrivial connected normal solvable subgroups. Equivalently, a reductive group is semisimple if and only if the identity component of its center is trivial.

\begin{exercise}
In this exercise we'll see some simple examples.
\begin{enumerate}[label=(\alph*)]
\item As a sanity check, give two examples of reductive groups that are not semisimple, two examples of semisimple groups with nontrivial center, and two examples of semisimple groups with trivial center.
\item Note that the additive group $\mathbb{C}^n$ is not reductive (because it is unipotent.) To see that its representations are generally not semisimple, find a representation of $\mathbb{C}^n$ that is not completely reducible.
\end{enumerate}
\end{exercise}

Exercise \ref{ex1} shows that any symmetric space $\xo$ is homogeneous for the action of a reductive group $H$. In the early 1980s, de Concini and Procesi \cite{dec.pro:83} observed that, if $H$ has trivial center, then $\xo$ has a distinguished equivariant compactification with remarkable properties.

\begin{definition}
\label{def:wond}
A \emph{wonderful compactification} of $\xo$ is an equivariant embedding $X$ such that
\begin{enumerate}[label=(\arabic*)]
\item $X$ is smooth and projective,
\item the boundary divisor
\begin{align*}
D&\coloneqq X\backslash\xo\\
&\,=D_1\cup\ldots\cup D_l
\end{align*}
has simple normal crossings and nonempty partial intersections, and
\item two points $a,b\in X$ are in the same $H$-orbit if and only if they are contained in exactly the same irreducible components of $D$---that is,
\[\{j\mid a\in D_j\}=\{j\mid b\in D_j\}.\]
\end{enumerate}
\end{definition}

\begin{theorem}\cite[Theorem 3.1]{dec.pro:83}
\label{decpro}
Suppose that the center of $H$ is trivial and that $\sigma:H\longrightarrow H$ is an involution. Then the symmetric space
\[H/H^\sigma\]
has a wonderful compactification.
\end{theorem}

\begin{exercise}
Suppose that $X$ is a wonderful variety. 
\begin{enumerate}[label=(\alph*)]
\item Show that, for any proper subset $J\subset\{1,\ldots,l\}$, the partial intersection
\[\bigcap_{j\not\in J}D_j\]
is the closure of a single $H$-orbit on $D$.
\item Conclude that $X$ contains a unique closed $H$-orbit of minimal dimension, given by 
\[\O_\varnothing=D_1\cap\ldots\cap D_l.\]
\end{enumerate}
\end{exercise}

Not all homogeneous spaces have compactifications that satisfy the properties of Definition \ref{def:wond}. For example, the next exercise shows that there is no such thing as a nontrivial ``wonderful toric variety.''

\begin{exercise}
\label{projaff}
Suppose that $X$ is a toric variety for the action of a complex torus $T$ which has a unique closed orbit. Show that $X$ is an affine variety. (Hint: Use the fact, due to Sumihiro, that any normal $T$-variety is covered by $T$-invariant affine open subsets \cite[Corollary 2]{sum:74}.)
\end{exercise}

It is therefore natural to ask which homogeneous spaces have wonderful compactifications. The answer to this question has its roots in a wide-ranging research program of Luna and Vust that began in the 1980s, which gave a systematic classification of a class of equivariant embeddings known as spherical varieties.

To formulate their definition, we will assume from now on that $H$ is a reductive group. A maximal solvable subgroup of $H$ is called a \emph{Borel subgroup}, and all such subgroups are conjugate.

\begin{example}
Suppose that $H=GL_n$. Then any Borel subgroup is conjugate to the subgroup
\[B=\left\{\begin{pmatrix}
			*& \ldots &* \\
			&\ddots&\vdots \\
			&&*
			\end{pmatrix}
			\right\}\]
of upper-triangular matrices. 
\end{example}

\begin{definition}
The homogeneous space $\xo$ is \emph{spherical} if a Borel subgroup of $H$ acts on $\xo$ with an open dense orbit. Equivariant embeddings of such spaces are called \emph{spherical varieties}.
\end{definition}

\begin{example}
\label{exs}
\begin{enumerate}[label=(\alph*), leftmargin=25pt, itemsep=5pt]
\item Since any torus $T$ is a solvable reductive group, all toric varieties are spherical.
\item 
\label{bruhat}
Let $H=G\times G$ for some reductive group $G$, and let $B\subset G$ be a Borel subgroup. Fix a maximal torus $T$ contained in $B$, let 
\[W=N_G(T)/T\]
be the associated Weyl group, and for each element $w\in W$ fix a coset representative $\dot{w}$ in $N_G(T)$. Then the Bruhat decomposition tells us that $G$ can be written as a disjoint union
\[G=\bigsqcup_{w\in W}B\dot{w}B\]
of finitely many $B\times B$-orbits, indexed by elements of $W$ and called \emph{Bruhat cells}. In particular, when $w_0\in W$ is the longest element of the Weyl group, the Bruhat cell
\[B\dot{w}_0B\]
is open and dense in $G.$ Since this is a single $B\times B$-orbit, it follows that $G$ itself is a spherical $G\times G$-space.
\item
\label{flag}
Let $G$ be a reductive group and let $B$ be a Borel subgroup. Then the projective homogeneous space 
\[G/B\]
is called the \emph{flag variety} of $G$. By the Bruhat decomposition of part \ref{bruhat}, it is spherical relative to the left action of $G$.
\end{enumerate}
\end{example}

Vust showed in \cite[Theorem 1]{vus:74} that all symmetric spaces are spherical. (The converse is not true---for example, the flag variety of Example \ref{exs}\ref{flag} is not a symmetric space.) Luna and Vust \cite{lun.vus:83} then developed a classification of equivariant embeddings of spherical spaces, using combinatorial objects known as \emph{colored fans}. 

Among the equivariant compactifications of a spherical homogeneous space $\xo$ there are two distinguished classes. The first is the class of \emph{toroidal compactifications}, which behave locally like toric varieties. Such compactifications correspond to colorless fans, and their boundaries are ``large,'' in the sense that the boundary divisor has simple normal crossings. Every equivariant compactification of $\xo$ is dominated by a toroidal one through a sequence blow-ups.

The second class is the class of \emph{simple compactifications}, which are those with a unique closed orbit. As we have seen in Exercise \ref{projaff}, not all homogeneous spaces admit simple compactifications. In fact, $\xo$ has a simple compactification if and only if the stabilizer $K$ of a point of $\xo$ has finite index in its normalizer \cite[Corollary 5.3]{bri.pau:87}. 

In this case, there is a unique simple compactification $X$ which is also toroidal. It dominates all other simple compactifications and is dominated by all toroidal ones. In other words, it has the universal property that, for every toroidal compactification $X'$ and every simple compactification $X''$ of $\xo$, there exist proper $G$-equivariant morphisms
\[X'\longrightarrow X\longrightarrow X''\]
which restrict to the identity along the open dense orbit $\xo$. If this distinguished compactification is smooth, then it is a wonderful variety.  Smoothness is guaranteed, for example, if the stabilizer $K$ is self-normalizing \cite[Corollary 7.2]{kno:96}, and this occurs in particular when $K$ is the fixed-point subgroup of an involution and $\xo$ is a symmetric space. \\

%
%
%
%
%
%
%
\section{The wonderful compactification of a semisimple algebraic group}
Let $G$ be a reductive group whose center is trivial---as a model, we will consider the running example of the projectivized general linear group $PGL_n$. We will describe one way to construct the wonderful compactification of $G$, viewed as a spherical homogeneous space for the action of the reductive group $G\times G$. This construction is an adaptation of the original one given by de Concini and Procesi \cite{dec.pro:83}, and is explained in detail in the now-classical survey \cite{eve.jon:08} by Evens and Jones.

Fix a Borel subgroup $B$ of $G$ containing a maximal torus $T$. In the simply-connected cover $\Gs$ of $G$, we write $\Ts$ and $\Bs$ for the subgroups which are their preimages. The irreducible representations of $\Gs$ can then be indexed by weights of the torus $\Ts$ which are dominant with respect to the Borel $\Bs$, and a representation is called \emph{regular} if its dominant weight does not lie on any walls of the corresponding Weyl chamber.

\begin{exercise}
\label{reg}
Let $\lambda$ be a dominant weight of $\Ts$ and let $V$ be the corresponding irreducible representation. If $v_\lambda\in V$ is a highest weight vector, show that the following are equivalent:
\begin{enumerate}[label=(\alph*)]
\item $\lambda$ is regular.
\item The stabilizer of the subspace $\mathbb{C}v_\lambda$ in $G$ is $B$.
\item The stabilizer of $\lambda$ in $W$ is trivial.
\end{enumerate}
\end{exercise}

If $V$ is a regular irreducible representation of $\Gs$, we can construct the following diagram:
\begin{equation*}
\begin{tikzcd}[row sep=large, column sep=large]
		\Gs	\arrow[r, "\rho"]\arrow{d}		& \text{End }V\backslash\{0\} \arrow{d}\\
		G 			\arrow[r, blue, hook, dotted, "\psi"]	& \mathbb{P}(\text{End }V).
\end{tikzcd}
\end{equation*}
The top arrow is the representation map, the left vertical arrow is the quotient by the center of $\Gs$, and the right vertical arrow is the quotient by scalars.

\begin{exercise}
Use Exercise \ref{reg} to show that the representation map $\rho$ descends to a well-defined morphism 
\[\psi:G\hooklongrightarrow \mathbb{P}(\text{End }V),\]
and that this morphism is an embedding.
\end{exercise}

\begin{definition}
The \emph{wonderful compactification} of $G$ is the closure
\[\Gbar\coloneqq \overline{\psi(G)}\]
of the image of the embedding $\psi$.
\end{definition}

The variety $\Gbar$ is a smooth projective variety which is independent of the choice of regular representation $V$. It is also stable under the $G\times G$-action on $\mathbb{P}(\text{End }V)$ given by pre- and post-composition of projectivized endomorphisms, which extends the action of $G$ on itself by left- and right-multiplication---therefore, it is an equivariant compactification in the sense of Definition \ref{def:equiv}.

The boundary divisor 
\begin{align*}
D&\coloneqq \Gbar\backslash G\\
	&=D_1\cup\ldots\cup D_l
	\end{align*}
is a union of smooth irreducible components whose number is equal to the rank $l$ of $G$. They have simple normal crossings and nonempty partial intersections---in other words, $\Gbar$ is a wonderful compactification of $G$.

The $G\times G$-orbits on $D$ are indexed by proper subsets of the simple roots $\{1,\ldots,l\}$. More precisely, to each proper subset $J\subset\{1,\ldots,l\}$ we can associate the data of a parabolic subgroup $P_J$, generated by the positive Borel $B$ and the negative simple roots 
\[\{\alpha_j\mid j\in J\}.\]
It has an opposite parabolic $P_J^-$, and we write $U_J^\pm$ for their unipotent radicals and $L_J\coloneqq P_J\cap P_J^-$ for their common Levi component. The Levi decomposition then gives
\[P_J^\pm=U_J^\pm L_J.\]

\begin{example}
\label{ex:para}
Let $G=PGL_4$ be the projectivized general linear group of rank $3$, and let $J=\{1,3\}$. The corresponding parabolic data is
\begin{align*}
&P_J=\left\{\begin{pmatrix}
			*& * & * &* \\
			*& * & * &* \\
			0& 0 & * &* \\
			0& 0 & * &* \\
			\end{pmatrix}
			\right\} \qquad
			P^-_J=\left\{\begin{pmatrix}
								*& * & 0 &0 \\
								*& * & 0 &0 \\
								*& * & * &* \\
								*& * & * &* \\
			\end{pmatrix}
			\right\}\\
&U_J=\left\{\begin{pmatrix}
			0& 0 & * &* \\
			0& 0 & * &* \\
			0& 0 & 0 &0 \\
			0& 0 & 0 &0 \\
			\end{pmatrix}
			\right\} \qquad
			U^-_J=\left\{\begin{pmatrix}
								0& 0 & 0 &0 \\
								0& 0 & 0 &0 \\
								*& * & 0 &0 \\
								*& * & 0 &0 \\
			\end{pmatrix}
			\right\}\\
&\hspace{1in} L_J=\left\{\begin{pmatrix}
								*& * & 0 &0 \\
								*& * & 0 &0 \\
								0& 0 & * &* \\
								0& 0 & * &* \\
			\end{pmatrix}
			\right\}
\end{align*}		
\end{example}

The index set $J$ corresponds to a $G\times G$-orbit $\O_J$ whose closure is the partial intersection
\[\overline{\O_J}=\bigcap_{j\not\in J}D_j\]
of divisor components not indexed by $J$. This orbit contains a distinguished basepoint $z_J$ whose stabilizer is 
\begin{equation}
\label{stabs}
\stab_{G\times G}(z_J)=\left\{(ul_1,vl_2)\mid u\in U_J,v\in U_J^-, l_1l_2^{-1}\in Z_{L_J}\right\}.
\end{equation}

\begin{exercise}
Show that the orbit $\O_J$ is a fibration
\begin{equation}
\label{fib1}
\begin{tikzcd}[row sep=large, column sep=large]
		G_J	\arrow[r, hook]	& \O_J \arrow[d]\\
								& G/P_J\times G/P_J^-
\end{tikzcd}
\end{equation}
over the product of partial flag varieties $G/P_J\times G/P_J^-$, with fiber is isomorphic to the quotient $G_J\coloneqq L_J/Z_{L_J}$. 
\end{exercise}

The group $G_J$ is itself a reductive group with trivial center, and when we take the closure of the orbit $\O_J$ we obtain a fibration
\begin{equation}
\label{fib2}
\begin{tikzcd}[row sep=large, column sep=large]
		\overline{G_J}	\arrow[r, hook]	& \overline{\O_J} \arrow[d]\\
								& G/P_J\times G/P_J^-,
\end{tikzcd}
\end{equation}
whose fiber is the wonderful compactification of $G_J$. In particular, when $J=\varnothing$ is the empty set, we obtain the unique closed $G\times G$-orbit
\[\O_\varnothing=D_1\cap\ldots\cap D_l.\]
In this case $P_\varnothing=B$ is a Borel subgroup, $P_\varnothing^-=B^-$ is the opposite Borel, and $L_\varnothing=T$ is the maximal torus of $G$, which is abelian. Therefore the subgroup $G_\varnothing$ is trivial, and diagrams \eqref{fib1} and \eqref{fib2} imply that this closed orbit
\[\O_\varnothing=\overline{\O_\varnothing}\cong G/B\times G/B^-\]
is isomorphic to the product of two copies of the full flag variety of $G$.

\begin{example}
\label{sl2}
Let $G=PGL_2$. Then $\Gs=SL_2$, and the weights of the maximal torus of are given by the integers. The irreducible representations are indexed by the nonnegative integers, and all nonzero weights are regular. Therefore we can take $V=\mathbb{C}^2$ to be the standard representation. In this case
\[\psi:PGL_2\hooklongrightarrow\mathbb{P}(M_{2\times 2})\]
is the standard embedding with image
\[\psi(PGL_2)=\left\{\begin{bmatrix} 
					a & b \\
					c & d
					\end{bmatrix}\in\mathbb{P}(M_{2\times 2})\mid ad-bc\neq0\right\}.\]
The wonderful compactification of $PGL_2$ is its closure
\[X=\mathbb{P}(M_{2\times2})\cong\mathbb{P}^3.\]

The boundary of this compactification is a smooth irreducible divisor with a transitive action of $PGL_2\times PGL_2$. It is given by the zero locus 
\[D=\left\{\begin{bmatrix} 
					a & b \\
					c & d
					\end{bmatrix}\in\mathbb{P}(M_{2\times 2})\mid ad-bc=0\right\}\]
of the determinant polynomial, and consists exactly of the projectivized $2$-by-$2$ matrices of rank $1$. In Exercise \ref{segre}, we will show that $D$ is isomorphic to a product of two copies of $\mathbb{P}^1$, which is the flag variety of $PGL_2$.
\end{example}

\begin{exercise}
\label{segre}
Construct an isomorphism
\[\mathbb{P}^1\times\mathbb{P}^1\xlongrightarrow{\sim} \left\{\begin{bmatrix} 
					a & b \\
					c & d
					\end{bmatrix}\in\mathbb{P}(M_{2\times 2})\mid ad-bc=0\right\}.\]
(You answer is probably an example of the \emph{Segre embedding}.)
\end{exercise}

Example \ref{sl2} does not generalize to larger ranks---when $n\geq 3$, the standard representation of $SL_n$ is not regular, and the wonderful compactification of $PGL_n$ is not simply the projective space $\mathbb{P}^{n^2-1}$. While this projective space is still a simple equivariant compactification of $PGL_n$, its boundary divisor is the vanishing locus 
\[\left\{A\in \mathbb{P}(M_{n\times n})\mid \det(A)=0\right\}\]
of the determinant polynomial, which is irreducible and has singularities along the subvarieties of $\mathbb{P}(M_{n\times n})$ where the rank of the matrix $A$ drops. 

\begin{exercise}
Check by hand that the compactification $\mathbb{P}(M_{3\times 3})$ of $PGL_3$ is not a wonderful variety.\\
\end{exercise}

%
%
%
%
%
%
%
\section{Logarithmic geometry of $\Gbar$}
Let $X$ be a smooth algebraic variety containing a simple normal crossings divisor $D$, and write $\xo$ for the complement of $D$ in $X$. The \emph{logarithmic tangent sheaf} or \emph{log-tangent sheaf} of $X$ with respect to $D$ is the subsheaf
\[\T_{X,D}\subset\T_X\]
consisting of vector fields on $X$ that are tangent to $D$. In other words, the sections of $\T_{X,D}$ are those derivations which preserve the ideal sheaf of local functions that vanish along $D$. 

Since $D$ has simple normal crossings, the sheaf $\T_{X,D}$ is locally free---if $x_1,\ldots,x_n$ are local coordinates on $X$ such that
\[D=\left\{(x_1,\ldots,x_n)\mid x_1\ldots x_k=0\right\},\]
the sheaf $\T_{X,D}$ is locally generated by the vector fields
\[x_1\partial_1,\ldots,x_k\partial_k, \partial_{k+1},\ldots,\partial_n.\]
This implies that $\T_{X,D}$ is the sheaf of sections of a vector bundle $T_{X,D}$, called the \emph{log-tangent bundle} of the pair $(X,D)$. Since the sheaves $\T_{X,D}$ and $\T_X$ have the same sections over the open dense locus $\xo$, the restriction of the bundle $T_{X,D}$ to $\xo$ is precisely the usual tangent bundle $T_{\xo}$.

\begin{exercise}
The inclusion $\T_{X,D}\hooklongrightarrow \T_X$ induces a morphism of vector bundles
\begin{equation}
\label{anchor}
T_{X,D}\longrightarrow T_X.
\end{equation}
Is this morphism injective? Surjective? Describe its image.
\end{exercise}

The dual vector bundle of $T_{X,D}$ is the \emph{log-cotangent bundle} $T^*_{X,D}$, whose sections are meromorphic forms with logarithmic poles along the divisor $D$. In the same coordinates, it is locally generated by the forms
\[\frac{dx_1}{x_1},\ldots,\frac{dx_k}{x_k}, dx_{k+1},\ldots,dx_n.\]
In particular, its restriction to $\xo$ is the cotangent bundle $T_{\xo}^*$.

The usual differential gives the sheaves 
\[\wedge^\bullet\T_{X,D}^*\]
the structure of a complex, and a global section $\omega\in\Gamma(\wedge^2T^*_{X,D})$ is called \emph{log-symplectic} if it is closed and logarithmically nondegenerate---that is, if
\begin{itemize}[topsep=3pt, itemsep=3pt]
\item $d\omega=0$, and
\item the induced map $\omega^\flat:T_{X,D}\longrightarrow T_{X,D}^*$ is an isomorphism.
\end{itemize}
In the next exercise, we will show that any log-cotangent bundle has a canonical log-symplectic structure.

\begin{exercise}
\label{logcan}
Let $M$ be a smooth algebraic variety with a simple normal crossings divisor $Z$, and write $\mathring{M}$ for the complement of $Z$ in $M$. Let $\pi:T^*_{M,Z}\longrightarrow M$ be the bundle map, and write $D\coloneqq\pi^{-1}(Z)$ for the preimage of the boundary divisor.
\begin{enumerate}[label=(\alph*)]
\item Check that $D$ is a simple normal crossings divisor in $T^*_{M,Z}$.
\item Show that the canonical symplectic structure on $T^*_{\mathring{M}}$ extends to a log-symplectic structure on the log-cotangent bundle $T^*_{M,Z}$. (Hint: Try to write a ``logarithmic'' version of the usual construction of the canonical symplectic structure on a cotangent bundle.)
\end{enumerate}
\end{exercise}

Once again let $H$ be a connected complex algebraic group that acts on $X$, and write $\h$ for its Lie algebra. The infinitesimal action of $\h$ gives a morphism of vector bundles
\begin{equation}
\label{action}
\rho:X\times\h\longrightarrow T_X
\end{equation}
from the trivial $\h$-bundle on $X$ to the tangent bundle $T_X$, which associates to each Lie algebra element $\xi$ the action vector field $\xi_X$. Writing $\O_X$ for the structure sheaf of $X$, on the level of sheaves the map $\rho$ corresponds to a morphism
\[\O_X\otimes\h\longrightarrow \T_X\]
that takes values in the tangent sheaf $\T_X$.

\begin{exercise}
\label{loghom}
Show that $X$ is homogeneous for the action of $H$ if and only if the bundle map \eqref{action} is surjective.
\end{exercise}

If the simple normal crossings divisor $D$ is stable under the action of $H$, the infinitesimal action map \eqref{action} factors through the morphism \eqref{anchor} via
\[\rho_D:X\times\h\longrightarrow T_{X,D}.\]
In other words, there is a commutative diagram
\begin{equation*}
\begin{tikzcd}[row sep=large, column sep=large]
		X\times\h	\arrow[r, "\rho_D"]\arrow[rd, swap, "\rho"]	& T_{X,D} \arrow[d]\\
								& T_X.
\end{tikzcd}
\end{equation*}

\begin{exercise}
Returning to the notation of Exercise \ref{logcan}, suppose that $H$ acts on $M$ and preserves $Z$. Then the induced action on the cotangent bundle $T^*_{\mathring{M}}$ is Hamiltonian, and for any covector $\alpha$ on $\mathring{M}$ the moment map
\[\mu:T^*_{\mathring{M}}\longrightarrow\h^*\]
is given by
\[\qquad\mu(\alpha)(\xi)=\alpha(\xi_M)\quad\text{for all }\xi\in\h.\]
Show that this moment map extends to a moment map
\[\mubar:T^*_{M,Z}\longrightarrow\h^*\]
with respect to the canonical log-symplectic structure on $T^*_{M,Z}$.
\end{exercise}

The pair $(X,D)$ is called \emph{log-homogeneous} if the morphism $\rho_D$ is surjective. In this case, by Exercise \ref{loghom}, the open dense locus $\xo$ is homogeneous for the action of $H$, and there is a short exact sequence of vector bundles 
\[0\longrightarrow\ker\rho_D\longrightarrow X\times\h\xlongrightarrow{\rho_D} T_{X, D}\longrightarrow 0.\]
Log-homogeneous varieties were introduced and studied by Brion \cite{bri:07}, who gave the following characterization of the kernel of the logarithmic action map $\rho_D$.

\begin{proposition}
\cite[Proposition 2.1.2]{bri:07}
Suppose that $(X,D)$ is log-homogeneous for the action of $H$, and fix a point $x\in X$ with $H$-orbit $\O=H\cdot x$. The kernel of $\rho_{D,x}$ is the Lie algebra of the kernel of the action of $\stab_H(x)$ on the normal space
\[N_{O,x}=T_{X,x}/T_{\O,x}.\]
\end{proposition}

All wonderful varieties are log-homogeneous \cite[Proposition 4.2]{pez:18}. In the special case when $H=G\times G$ and $X=\Gbar$, the logarithmic action map gives a short exact sequence of vector bundles
\[0\longrightarrow\ker\rho_D\longrightarrow\Gbar\times\g\times\g\xlongrightarrow{\rho_D} T_{\Gbar, D}\longrightarrow 0.\]

\begin{exercise}
\label{kernels}
\begin{enumerate}[label=(\alph*)]
\item Let $e\in G$ be the identity element. Show that the kernel of $\rho_D$ at $e$ is the diagonal subalgebra
\[\g_\Delta=\left\{(x,x)\in\g\times\g\mid x\in\g\right\}.\]
\item Let $a\in G$ be an arbitrary group element. Find the kernel of $\rho_D$ at $a$.
\item Let $z_J$ be the distinguished basepoint of the $G\times G$-orbit $\O_J$, with stabilizer given by \eqref{stabs}. Show that the kernel of $\rho_D$ above $z_J$ is the subalgebra
\[\p_J\times_{\l_J}\p_J^-,\]
where $\p_J^\pm$ denotes the Lie algebra of the parabolic subgroup $P_J$ and $\l_J$ is the Lie algebra of the common Levi $L_J$.
\end{enumerate}
\end{exercise}

Now we consider the Killing form
\begin{align*}
\kappa:\g\times\g&\longrightarrow\mathbb{C}\\
		(x,y)&\longmapsto \text{tr}(\ad_x\circ\ad_y),
		\end{align*}
and we use it to define a symmetric, nondegenerate, invariant bilinear form on $\g\times\g$ by
\begin{equation}
\label{inner}
\langle(x_1,x_2),(y_1,y_2)\rangle=\kappa(x_1,y_1)-\kappa(x_2,y_2).
\end{equation}
With respect to this form, the diagonal subalgebra $\g_\Delta$ is Lagrangian. Since $\kappa$ is $G$-invariant, this implies that the subalgebra
\[\ker\rho_{D,a}=(a,e)\cdot\g_\Delta\]
of Exercise \ref{kernels} is also Lagrangian for any $a\in G$. Moreover, since the condition of being Lagrangian is closed, this means that the bundle
\[\ker\rho_D\subset\Gbar\times\g\times\g\]
is a Lagrangian subbundle. 

We therefore obtain a natural isomorphism between $\ker\rho_D$ and the dual of the quotient
\[\left((\Gbar\times\g\times\g)/\ker\rho_D\right)^*\cong T^*_{\Gbar,D}.\]
In other words, the bundle $\ker\rho_D$ is isomorphic to the log-cotangent bundle of the wonderful compactification $\Gbar$. Our short exact sequence becomes
\begin{equation}
\label{les}
0\longrightarrow T^*_{\Gbar,D}\longrightarrow\Gbar\times\g\times\g\longrightarrow T_{\Gbar, D}\longrightarrow 0,
\end{equation}
and this realizes the log-cotangent bundle as a bundle of Lie algebras over $\Gbar$. Note that this vector bundle is not locally trivial as a bundle of Lie algebras---the fiber above the identity element is isomorphic to the semisimple Lie algebra $\g_\Delta,$ while the fibers over the $G\times G$-orbit $\O_J$ are not semisimple, as the next exercise shows.

\begin{exercise}
Check that the Lie algebra $\p_J\times_{\l_J}\p_J^-$ is not reductive. (Hint: Reductive Lie algebras contain no nontrivial nilpotent ideals.)
\end{exercise}

Restricting the sequence \eqref{les} to the open dense locus $G$ of $\Gbar$, we obtain the commutative diagram
\begin{equation*}
\begin{tikzcd}[row sep=large, column sep=large]
		0\arrow[r] 	&T^*_{\Gbar,D}	\arrow[r, hook]	& \Gbar\times\g\times\g\arrow[r] &T_{\Gbar,D}\arrow[r]	&0\\
				&T^*_G\arrow[r, hook]\arrow[u] & G\times\g\times\g,\arrow[u] &&
\end{tikzcd}
\end{equation*}
in which the square is Cartesian. If we identify the cotangent bundle $T^*_G$ with $G\times\g$ using left-trivialization and the isomorphism between $\g$ and $\g^*$ given by the Killing form, the bottom inclusion is given by the map
\begin{align*}
T^*_G\cong G\times \g&\hooklongrightarrow G\times\g\times\g\\
			(a,x)&\longmapsto (a,\Ad_ax,x).
			\end{align*}
The last two components of this morphism recover the moment map
\begin{align*}
\mu:T^*_G\cong G\times\g&\longrightarrow \g\times\g\\
				(a,x)&\longmapsto (\Ad_ax,x)
				\end{align*}
for the Hamiltonian $G\times G$-action on $T^*G$ given by
\[(g,h)\cdot(a,x)=(gah^{-1}, \Ad_hx).\]
Projection onto the fibers of the trivial bundle $\Gbar\times\g\times\g$ then extends this to the moment map
\[\overline{\mu}:T^*_{\Gbar,D}\longrightarrow\g\times\g\]
relative to the canonical log-symplectic structure on $T^*_{\Gbar,D}$ constructed in Exercise \ref{logcan}.

Finally, we note that the short exact sequence \eqref{les} induces a morphism
\begin{align}
\label{demazure}
\Gbar&\longrightarrow \gr(\dim\g,\g\oplus\g)\\
	a&\longmapsto [\ker\rho_{D,a}],\nonumber
	\end{align}
from the wonderful compactification $\Gbar$ into the Grassmannian of half-dimensional subspaces of $\g\oplus\g$. This morphism is an embedding, called the \emph{Demazure embedding} \cite{los:09}. By what we saw above, its image lies in the closed subvariety 
\[\L(\g\oplus\g)\subset\gr(\dim\g,\g\oplus\g)\]
consisting of subalgebras of $\g\oplus\g$ which are Lagrangian with respect to the inner product \eqref{inner}. This subvariety has many remarkable properties, which we now study in a slightly more general context.\\

%
%
%
%
%
%
%
\section{Lagrangian subalgebras and Poisson-homogeneous spaces}
Let $\d$ be a real or complex Lie algebra of dimension $2n$ that carries a nondegenerate, symmetric, invariant bilinear form $\langle\cdot,\cdot\rangle.$ If
\[\d=\l_1+\l_2\]
is a splitting of $\d$ into Lagrangian subalgebras and we write $\uu=\l_1$, there is a natural identification of $\uu^*$ with $\l_2$. Then
\[(\d,\uu,\uu^*)\]
is a Manin triple, and the Lie bracket on $\uu^*$ gives a Poisson bracket $\pi_U$ on any Lie group $U$ integrating $\uu$. This bracket has the property that the group multiplication map is a Poisson morphism---in other words, $(U,\pi_U)$ is a \emph{Poisson--Lie group}. Equivalently, the Poisson bivector $\pi_U$ satisfies the identity
\[\pi_{U,gh}=L_g\pi_{U,h}+R_h\pi_{U,g}\quad\text{for all }g,h\in U,\]
where $L_g$ and $R_h$ are the differentials of left multiplication by $g$ and right multiplication by $h$, respectively.

\begin{example}
\label{reals}
Let $\d=\g$ be a semisimple complex Lie algebra of complex dimension $n$, viewed as a $2n$-dimensional Lie algebra over the real numbers. Consider the nondegenerate, symmetric, invariant bilinear form 
\[\langle\cdot,\cdot\rangle=\text{Im}(\kappa)\]
given by the imaginary part of the Killing form. The Iwasawa decomposition
\[\g=\mathfrak{k}+\mathfrak{a}+\mathfrak{n}\]
gives a Lagrangian splitting of $\g$ into the Lagrangian subalgebras $\uu=\mathfrak{k}$ and $\uu^*=\mathfrak{a}+\mathfrak{n}$, and 
\[(\g,\k,\a+\n)\]
is a Manin triple. In this way, any real Lie group $K$ integrating $\k$ becomes a Poisson--Lie group. 
\end{example}

\begin{example}
\label{gbar}
Let $\g$ be any real or complex Lie algebra, and view $\g^*$ as an abelian Lie algebra equipped with the trivial Lie bracket. The semidirect product $\d=\g\ltimes\g^*$ with respect to the coadjoint action is the Lie algebra with bracket given by
\[[(x,\alpha),(y,\beta)]=([x,y],\ad_{x}\beta-\ad_{y}\alpha),\qquad\text{for all }(x,\alpha),(y, \beta)\in\g\ltimes\g^*.\]
It has a natural nondegenerate symmetric invariant inner product
\[\langle(x,\alpha),(y,\beta)\rangle=\alpha(y)+\beta(x)\qquad\text{for all }(x,\alpha),(y,\beta)\in\g\ltimes\g^*,\]
with respect to which the subalgebras $\g\oplus 0$ and $0\oplus\g^*$ are Lagrangian. If $G$ is a Lie group integrating $\g$, the double $\g\ltimes \g^*$ integrates to the group $G\ltimes\g^*$ with multiplication given by
\[(g,\alpha)\cdot(h,\beta)=(gh,\Ad_{h^{-1}}\alpha+\beta)\qquad\text{for all }(g,\alpha),(h,\beta)\in G\ltimes\g^*.\]
The Manin triple
\[(\g\ltimes\g^*,\g,\g^*)\]
induces a trivial Poisson structure on the group $G$, and the Kostant--Kirillov--Souriau Poisson structure on its Poisson--Lie dual group $\g^*$.
\end{example}

\begin{example}
\label{evlu2}
Let $\g$ be a complex semisimple Lie algebra and fix a maximal torus $\t$ and a Borel subalgebra $\b$ that contains it. Writing $\b^-$ for the opposite Borel we obtain a Lagrangian splitting
\begin{equation}
\label{stdsplit}
\g\oplus\g=\g_\Delta+\b\times_{\t}\b^-
\end{equation}
with respect to the nondegenerate invariant bilinear form
\[\langle(x_1,x_2),(y_1,y_2)\rangle=\kappa(x_1,y_1)-\kappa(x_2,y_2)\]
defined in \eqref{inner}. This produces a Manin triple
\[(\g\oplus\g,\g_{\Delta},\b\times_{\t}\b^-),\]
and gives a canonical Poisson--Lie group structure on the group $G$ integrating $\g$. 
\end{example}

Poisson--Lie groups were introduced by Drinfeld \cite{dri:83} and by Semenov-Tian-Shansky \cite{sts:85}. For background, we refer to \cite[Section 1]{cha.pre:95} or to \cite[Section 1]{lu.wei:90}. Given such a group $(U,\pi_U)$, it is natural to consider its \emph{Poisson actions}---that is, actions of $U$ on a Poisson manifold $(M,\pi)$ that have the property that the action map
\[U\times M\longrightarrow M\]
is a Poisson map. Equivalently, this means that the Poisson bivector $\pi$ satisfies the identity
\[\pi_{gm}=g_*\pi_{m}+m_*\pi_{U,g}\quad\text{for all }g\in U, m\in M,\]
where $g_*$ and $m_*$ are the differentials of the maps
\begin{align*}
g:M&\longrightarrow M\hspace{.7in}\text{and}\hspace{-.7in}& m:U&\longrightarrow M \\
m'&\longmapsto gm'& 		g'&\longrightarrow g'm
\end{align*}
respectively.

\begin{definition}
A \emph{$(U,\pi_U)$-homogeneous space} is a Poisson manifold $(M,\pi)$ equipped with an action of $U$ that is transitive and Poisson.
\end{definition}

Let $(M,\pi)$ be a $(U,\pi_U)$-homogeneous space, and for any point $m$ of $M$ write $U_m$ for the stabilizer of $m$ and $\uu_m$ for its Lie algebra. Then there is a natural identitication
\[T_{M,m}\cong\uu/\uu_m,\]
and we can view the Poisson bivector as $\pi_m$ as an element of $\wedge^2(\uu/\uu_m).$ This induces a natural map
\[\pi_m^\#:(\uu/\uu_m)^*\longrightarrow (\uu/\uu_m),\]
and we can define the subspace
\[\l_ m\coloneqq\left\{(x,\xi)\in \uu+\uu^*\mid \xi_{\vert\uu_m}=0\text{ and }\pi_m^\#(\xi)=x+\uu_m\right\}\]
 of the Lie algebra $\d$. In other words, $\l_m$ consists of pairs $(x,\xi)$ with the property that $\xi$ descends to an element of $(\uu/\uu_m)^*$ whose image under $\pi_m^\#$ is the coset of $x$. 
 
\begin{exercise}\cite[Theorem 1]{dri:93}
\label{lags}
Show that, for any point $m\in M$, the subspace $\l_m$ is a Lagrangian subalgebra of $\d$.
\end{exercise}
 
In this way we have obtained a map
\begin{align}
\label{lagmap}
\Psi:M&\longrightarrow \gr(n,\d)\\
m&\longmapsto [\l_m],\nonumber
\end{align}
known as the \emph{Drinfeld map} \cite{dri:93}, from the Poisson-homogeneous space $(M,\pi)$ to the Grassmannian of half-dimensional subspaces of the Lie algebra $\d$. Exercise \ref{lags} implies that its image is contained in the \emph{variety of Lagrangian subalgebras} 
\[\L(\d)\coloneqq\left\{[\l]\in\gr(n,\d)\mid [\l,\l]\subset\l\text{ and }\l^\perp=\l\right\}.\]
Since both of the defining conditions on the right are algebraic, $\L(\d)$ is a closed subvariety of the Grassmannian. Moreover, it is possibly singular, possibly disconnected, and stabilized by the adjoint action of any algebraic group $D$ that integrates $\d$.
 
\begin{theorem}\cite[Theorem 1]{dri:93}
\label{drinfeld}
Let $(M,\pi)$ be a $(U,\pi_U)$-homogeneous space. 
\begin{enumerate}[label=\textup{(\alph*)}]
\item The Drinfeld map $\Psi$ is $U$-equivariant---in other words, for any $m\in M$ and $u\in U$, 
\[\l_{um}=\Ad_u(\l_m).\]
\item The assignment $(M,\pi)\longmapsto U\cdot[\l_m]$ defines an injection
\[\quad\qquad\left\{\text{$(U,\pi_U)$-homogeneous spaces with connected stabilizers}\}\hooklongrightarrow\{\text{$U$-orbits on $\L(\d)$}\right\}.\]
Its image consists of those Lagrangian subalgebras $\l$ of $\d$ with the property that $\l\cap\uu$ integrates to a closed subgroup of $U$.
\end{enumerate}
\end{theorem}

The variety $\L(\d)$ has a natural Poisson structure which is compatible with the map \eqref{lagmap}, and which was first introduced and studied by Evens and Lu \cite{eve.lu:01,eve.lu:06}. Before defining it we remark that, since $\L(\d)$ is a possibly singular variety, by a Poisson structure we mean a Lie bracket on the structure sheaf of $\L(\d)$ which is a derivative in each component. On the smooth locus of $\L(\d)$, this definition agrees with the usual notion of a Poisson bivector field.

Identify
\[\d^*\cong\uu^*+\uu\]
and define a classical $r$-matrix $R\in\wedge^2\d$ by
\[R((\xi_1,x_1),(\xi_2,x_2))=\xi_2(x_1)+\xi_1(x_2)\quad\text{for all }(\xi_1,x_i)\in\d^*.\]
If $D$ is an algebraic group integrating $\d$, the adjoint action of $D$ on $\gr(n,\d)$ induces an infinitesimal action map
\[\rho:\d\longrightarrow\Gamma(T_{\gr(n,\d)}).\]
Via this action, the $r$-matrix $R$ becomes a bivector field on $\gr(n,\d)$ given by
\[\Pi\coloneqq\rho(R).\]

\begin{exercise}
\label{jacobi}
Consider the identification 
\begin{align*}
\#:\quad\d^*\,\,\,&\longrightarrow\,\,\,\,\d\\
	(\xi,x)&\longmapsto(x,\xi).
	\end{align*}
\begin{enumerate}[label=(\alph*)]
\item Show that the Schouten bracket $[R,R]\in\wedge^3\d$ is given by
\[[R,R](d_1,d_2,d_3)=\langle\#d_1,[\#d_2,\#d_3]\rangle\quad\text{for all }d_i\in\d^*.\]
\item Suppose that $\l$ is a Lagrangian subalgebra of $\d$. Use part (a) to show that
\[[\Pi,\Pi]_{[\l]}=0.\]
\end{enumerate}
\end{exercise}

Since $\L(\d)$ is stable under the action of $D$, the skew-symmetric bracket induced by $\Pi$ on the structure sheaf of $\gr(n,\d)$ descends to a bracket on $\L(\d)$. Exercise \ref{jacobi} implies that this bracket satisfies the Jacobi identity, and therefore $\L(\d)$ has the structure of a Poisson algebraic variety. The $D$-orbits of $\L(\d)$ are Poisson submanifolds with respect to this structure.

\begin{remark}
Note that the Poisson structure $\Pi$ depends on the fixed Lagrangian splitting. There are therefore many possible Poisson structures on $\L(\d)$, given by different choices of Lagrangian splitting, and it is an interesting problem to classify them. In the case when $\d=\g\oplus\g$, such a classification is carried out by Lu and Yakimov in \cite{lu.yak:08}.
\end{remark} 

\begin{theorem}\cite[Proposition 2.16 and Theorem 2.17]{eve.lu:01}
\label{thm1}
\begin{enumerate}[label=\textup{(\alph*)}]
\item The action of $(U,\pi_U)$ on $(\L(\d),\Pi)$ is a Poisson action.
\item Every $U$-orbit on $\L(\d)$ is a Poisson submanifold of $\L(\d)$, and therefore a $(U,\pi_U)$-homogeneous space.
\end{enumerate}
\end{theorem}

Let $\l$ be a Lagrangian subalgebra of $\d$. According to Theorem \ref{thm1}, the orbit 
\[M=U\cdot[\l]\subset\gr(n,\d)\]
is a $(U,\pi_U)$-homogeneous space. The image of its basepoint under the Drinfeld map \eqref{lagmap} is the Lagrangian subalgebra
\[\Psi([\l])=[\uu_{[\l]}+(\uu+\uu_{[\l]}^\circ)\cap\l]\]
\cite[Theorem 2.20]{eve.lu:01}. In particular, the Drinfeld map in this case is not necessarily the identity map. In fact, viewed as a function from $\L(\d)$ to itself, the Drinfeld map is generally not even continuous. \cite[Remark 2.21]{eve.lu:01}.

\begin{exercise}
\begin{enumerate}[label=(\alph*)]
\item Show that $\Psi([\l])=[\l]$ if and only if 
\[\l\cap\uu=\uu_{[\l]}.\]
The Lagrangian subalgebras satisfying this property are called \emph{model points} of $\L(\d)$.
\item Check that the set of model points of $\L(\d)$ is stable under the action of $U$.
\item Suppose that $(M,\pi)$ is a $(U,\pi_U)$-homogeneous space and that, for some $m\in M$, $[\l_m]$ is a model point of $\L(\d)$. Show that the Drinfeld map 
\[\Psi:M\longrightarrow U\cdot[\l_m]\]
is a local diffeomorphism.
\end{enumerate}
\end{exercise}

\begin{theorem}\cite[Theorem 2.22]{eve.lu:01}
\label{drpoisson}
For every $(U,\pi_U)$-homogeneous space $(M,\pi)$, the Drinfeld map given by \eqref{lagmap} is a Poisson map.
\end{theorem}

Putting all of this together, Theorem \ref{drpoisson} shows that the Drinfeld map \eqref{lagmap} is a Poisson map from any given $(U,\pi_U)$-homogeneous space to a $U$-orbit in $\L(\d)$, which is itself a $(U,\pi_U)$-homogeneous space. When this homogeneous space consists of model points of $\L(\d)$, the Drinfeld map is a local Poisson diffeomorphism---in other words, model point $U$-orbits on the variety of Lagrangian subalgebras can be viewed as ``local models'' of Poisson homogeneous spaces.

\begin{exercise}
Consider the Manin triply $(\g,\k,\a+\n)$ of Example \ref{reals}, and let $\l=\a+\n$.
\begin{enumerate}[label=(\alph*)]
\item Show that $\l$ is not a model point of $\L(\g)$.
\item Show that the image of $\l$ under the Drinfeld map $\Psi$ is a model point of $\L(\d)$.\\
\end{enumerate}
\end{exercise}

%
%
%
%
%
%
%
\section{Poisson structures on wonderful spaces}
The variety of Lagrangian subalgebras $\L(\d)$ is an interesting object whose Poisson geometry has many unexpected connections to the representation theory of $\d$. Many interesting equivariant embeddings of homogeneous spaces appear as closures of $D$-orbits or $U$-orbits on $\L(\d)$. We conclude this survey with three illustrative examples, following once again the work of Evens and Lu \cite{eve.lu:01, eve.lu:06}. 

\begin{example}
Let $\g$ be a semisimple complex Lie algebra with the Lagrangian splitting given by the Iwasawa decomposition
\[\g=\mathfrak{k}+\mathfrak{a}+\mathfrak{n},\]
as in Example \ref{reals}. The variety of Lagrangian subalgebras $\L(\g)$ and its induced Poisson structure were studied by Evens and Lu in \cite{eve.lu:01}, using the parametrization of Lagrangian subalgebras of $\g$ given by Karolinsky \cite{kar:96}. 

In this work the authors classify the irreducible components of $\L(\g)$ and show that each one is a fibration over a partial flag variety with fiber the real points of the wonderful compactification associated to a Levi subgroup of $\g$. A large number of familiar $G$- and $K$-homogeneous spaces appear as Poisson submanifolds of $\L(\g)$. Among them are the flag varieties of $G$, equipped with the Bruhat--Poisson structure whose symplectic leaves are Bruhat cells \cite{lu.wei:90}.
\end{example}

\begin{example}
Let $\g$ be a real or complex Lie algebra, and consider the Manin triple
\[(\g\ltimes\g^*,\g,\g^*)\]
of Example \ref{gbar}. The normalizer of the Lagrangian subalgebra $\g^*$ in $G\ltimes\g^*$ is the subgroup
\[\stab_{G\ltimes\g^*}[\g^*]=G\ltimes\{0\},\]
and therefore the $G\ltimes\g^*$-orbit 
\[(G\ltimes\g^*)\cdot[\g^*]\cong (G\ltimes\g^*)/(G\times\{0\})\cong\g^*\]
is a Poisson submanifold of $\L(\g\ltimes\g^*)$. The induced Poisson structure on $\g^*$ is precisely the Kostant--Kirillov--Souriau Poisson structure.

The variety $\L(\g\ltimes\g^*)$ is a reducible algebraic variety whose irreducible components are often singular, and generally its geometry is very complicated. However, in the case when $\g$ is a complex semisimple Lie algebra, its points are indexed by a parametrization due to Karolinsky and Stolin \cite{kar.sto:02}. Recently, Evens and Li \cite{eve.li:20} used this parametrization to show that the closed $G\ltimes\g^*$-orbits of $\L(\g\ltimes\g^*)$ are in bijection with the abelian ideals of a maximal Borel of $\g$.
\end{example}

\begin{example}
Let $\g$ be a semisimple complex Lie algebra, and consider the Manin triple
\[(\g\oplus\g, \g_{\Delta}, \b\times_{\t}\b^-)\]
of Example \ref{evlu2}, which can be viewed as a complexification of Example \ref{reals}. This triple induces the standard Poisson--Lie group structure on the adjoint group $G$, and by the discussion in the previous section it gives the variety $\L(\g\oplus\g)$ a Poisson structure for which the orbits of $G\times G$, $G_\Delta$, and $B\times_T B^-$ are Poisson submanifolds. The geometry of this Poisson structure was studied by Evens and Lu in \cite{eve.lu:06}, and we devote the rest of this section to some of their results.
\end{example}

In \cite{eve.lu:06}, the authors classify the $G\times G$-orbits and the irreducible and connected components of $\L(\g\oplus\g)$ using generalized Belavin--Drinfeld triples, which were introduced by Schiffman \cite{sch:98} and which we now briefly recall.

\begin{definition}
A \emph{generalized Belavin--Drinfeld triple} is a triple $(I,J,\eta)$, where
\begin{itemize}
\item $I$ and $J$ are subsets of the set $\Delta$ of simple roots, and
\item $\eta:I\longrightarrow J$ is a bijection that preserves the inner product $\langle\cdot,\cdot\rangle$.
\end{itemize}
\end{definition}

The subsets $I$ and $J$ correspond to parabolic subgroups $P_I$ and $P_J$ as illustrated in Example \ref{ex:para}. We denote by $P_I^-$ and $P_J^-$ their opposite parabolics, and by $L_I$ and $L_J$ their respective Levi components. The isometry $\eta$ induces a group isomorphism between the adjoint forms $L_I/Z_{L_I}$ and $L_J/Z_{L_J}$, which we denote by $G_I$ and $G_J$ respectively.

\begin{theorem}\cite[Proposition 2.17 and Proposition 2.25]{eve.lu:06}
\label{orbsL}
Let $\O$ be a $G\times G$-orbit in $\L(\d)$.
\begin{enumerate}[label=\textup{(\alph*)}]
\item There is a generalized Belavin--Drinfeld triple $(I,J,\eta)$ such that $\O$ is a fibration 
\begin{equation*}
\begin{tikzcd}[row sep=large, column sep=large]
		G_I	\arrow[r, hook]	& \O \arrow[d]\\
								& G/P_I\times G/P_J^-
\end{tikzcd}
\end{equation*}
over the partial flag variety $G/P_I\times G/P_J^-$ with fiber isomorphic to the adjoint group $G_I$.
\item The closure $\overline{\O}$ of $\O$ is a fibration
\begin{equation*}
\begin{tikzcd}[row sep=large, column sep=large]
		\overline{G_I}	\arrow[r, hook]	& \overline{\O} \arrow[d]\\
								& G/P_I\times G/P_J^-
\end{tikzcd}
\end{equation*}
with fiber isomorphic to the wonderful compactification $\overline{G_I}$ of $G_I$.
\end{enumerate}
\end{theorem}

Theorem \ref{orbsL} implies that the closures of $G\times G$-orbits in $\L(\g\oplus\g)$ are smooth. Moreover, Evens and Lu show that each $G\times G$-orbit is a union of finitely many $B\times B^-$-orbits. Since one of these $B\times B^-$-orbits must be open and dense, this implies that every $G\times G$-orbit in $\L(\g\oplus\g)$ is a spherical $G\times G$-homogeneous space \cite[Corollary 2.22]{eve.lu:06}. 

\begin{example}
The normalizer of the Lagrangian subalgebra $\g_\Delta$ is the diagonal subgroup
\[\stab_{G\times G}[\g_\Delta]=G_\Delta.\]
This construction therefore gives a Poisson structure on the group $G$ itself, via the isomorphism
\[(G\times G)\cdot[\g_\Delta]\cong (G\times G)/G_\Delta\cong G,\]
which is different from the Poisson--Lie group structure on $G$. By the above discussion, this Poisson structure extends to the closure of the $G\times G$-orbit of $[\g_\Delta]$, which by Exercise \ref{kernels} is isomorphic to the wonderful compactification $\Gbar$ via the map \eqref{demazure}. Using this observation, the characterization of $G\times G$-orbits on $\Gbar$ given in \eqref{fib1} and \eqref{fib2} can also be deduced from the statement of Theorem \ref{orbsL}.
\end{example}

Now let $\O_\Delta$ be a $G_\Delta$-orbit, let $\O_{B\times B-}$ be a $B\times B^-$-orbit, and suppose that the two have nonempty intersection. If $\O$ is the $G\times G$-orbit that contains them both, the Lagrangian splitting \eqref{stdsplit} implies that $\O_\Delta$ and $\O_{B\times B^-}$ intersect transversally in $\O$. Therefore
\[\O_\Delta\cap\O_{B\times B^-}\]
is a smooth submanifold of $\L(\g\oplus\g)$ and, since it is an intersection of Poisson submanifolds, it is in fact a Poisson submanifold. 

\begin{theorem}\cite[Theorem 4.5 and Theorem 4.14]{eve.lu:06}
The intersection $\O_\Delta\cap\O_{B\times B^-}$ is a regular Poisson manifold---that is, the Poisson bivector has constant rank, so that the symplectic foliation is a regular foliation. The diagonal torus $T_\Delta$ acts transitively on the symplectic leaves of this intersection.
\end{theorem}

\begin{remark}
The intersections $\O_\Delta\cap\O_{B\times B^-}$ are called the \emph{$T$-leaves} of the Poisson structure on $\L(\g\oplus\g)$, because they are $T$-saturations of symplectic leaves. In the more general setting of complex tori acting on Poisson manifolds, such leaves have recently come to play an important role in cluster geometry \cite{lu.mou:17, lu.ele:21}.
\end{remark}

\begin{example}
In particular, when we view the group $G$ as a Poisson submanifold of $\L(\g\oplus\g)$, the $T$-leaves are exactly the intersections of the conjugacy classes with the Bruhat cells.
\end{example}

We conclude with an important class of $G_\Delta$-orbits on $\L(\g\oplus\g)$ which can be identified with the symmetric $G$-spaces defined in Definition \ref{symspaces}.

\begin{exercise}
\label{symlag}
Let
\[\sigma:\g\longrightarrow\g\]
be an algebraic involution of the Lie algebra $\g$, and consider its graph
\[\l_\sigma\coloneqq\{(x,\sigma(x))\in\g\oplus\g\mid x\in\g\}.\]
\begin{enumerate}[label=(\alph*)]
\item Show that $\l_\sigma$ is a Lagrangian subalgebra of $\g\oplus\g$.
\item Show that the stabilizer $\stab_{G_\Delta}(\l_\sigma)$ is isomorphic to the fixed-point subgroup $G^\sigma.$
\end{enumerate}
\end{exercise}

It follows from Exercise \ref{symlag} that the $G_\Delta$-orbit of the Lagrangian subalgebra $[\l_\sigma]$ in $\L(\g\oplus\g)$ is isomorphic to the symmetric space
\[G_\Delta\cdot[\l_{\sigma}]\cong G/G^\sigma.\]
Its closure in the variety of Lagrangian subalgebras is precisely the wonderful compactification of $G/G^\sigma$ given by Theorem \ref{decpro} \cite[Proposition 3.21]{eve.lu:06}. In this way, all de Concini--Procesi wonderful compactifications sit inside $\L(\g\oplus\g)$ as Poisson submanifolds, and their resulting Poisson structure is characterized by the broader results of Evens and Lu.

%
%
%
%
%
%
%
\bibliographystyle{plain}
\bibliography{biblio}

\end{document}